\documentclass[10pt]{article}

\usepackage{epsfig}
\usepackage{graphicx}
\usepackage{amsmath}
\usepackage{amssymb}
\usepackage{amsthm}

\title{On the relationship of continuity and boundary regularity in PMC Dirichlet problems}
\author{Kirk E. Lancaster and Jaron Melin                \\
Department of Mathematics, Statistics and Physics \\
                            Wichita State University \\
                            Wichita, Kansas, 67260-0033 \\
                            U.S.A.}

\def\Real{{\rm I\hspace{-0.2em}R}}
\def\Natural{{\rm I\hspace{-0.2em}N}}

\def\rest{
    \kern.2em\vrule width .06 em height 1.75 ex depth .20 ex
   \vrule height -.08 ex width .375 em depth .20 ex \kern .200 em}

\newcommand{\suchthat}{%
\mathrel{\ooalign{$\ni$\cr\kern-1pt$-$\kern-6.5pt$-$}}}

\newtheorem{thm}{Theorem}

\newtheorem{lem}{Lemma}

\begin{document}
\maketitle

\begin{abstract}
In 1976, Leon Simon showed that if a compact subset of the boundary of a domain is smooth and has negative mean curvature, then the non-parametric 
least area problem with Lipschitz continuous Dirichlet boundary data has a generalized solution which is continuous on the union of the domain and this compact 
subset of the boundary, even if the generalized solution does not take on the prescribed boundary data.  Simon's result has been extended to boundary value problems 
for prescribed mean curvature equations by other authors.   
In this note, we construct Dirichlet problems in domains with corners and demonstrate that the variational solutions of these Dirichlet problems are discontinuous 
at the corner, showing that Simon's assumption of regularity of the boundary of the domain is essential.  
\end{abstract}
\vspace{2mm}

\noindent keywords:  Nonconvex corner, prescribed mean curvature.

\noindent 2010 MSC: Primary 35J93, 76B45 (Secondary 35J62, 53A10)
\vspace{3mm}

\section{Introduction}

Let $n\in \Natural$  with $n\ge 2$  and  suppose $\Omega$  is a bounded, open set in $\Real^{n}$  with locally Lipschitz boundary $\partial\Omega.$   
Fix $H\in C^{2}\left(\Real^{n}\times\Real\right)$  such that $H$  is bounded and $H(x,t)$  is nondecreasing in $t$  for $x\in\Omega.$
Consider the prescribed mean curvature Dirichlet problem of finding a function $f\in C^{2}\left(\Omega\right)\cap  C^{0}\left(\overline{\Omega}\right)$  
which satisfies 
\begin{eqnarray}
\label{ONE-A}
{\rm div}\left( Tf\right) &  = &  H(x,f) \ \ \ \ \ {\rm in} \ \ \Omega, \\ 
f & = & \phi \ \ \ \ \ {\rm  on} \ \ \partial\Omega, 
\label{ONE-B}
\end{eqnarray}
where  $Tf= \frac{\nabla f}{\sqrt{1+\left|\nabla f\right|^{2}}}$  and $\phi\in C^{0}\left( \partial \Omega\right)$  is a prescribed function; 
such a function $f,$  if it exists, is a classical solution of the Dirichlet problem.  
It has been long known (e.g. Bernstein in 1912) that some type of boundary curvature condition (which depends on $H$)  must be satisfied 
in order to guarantee that a classical solution exists for each $\phi\in C^{0}\left( \partial \Omega\right)$  
(e.g. \cite{JenkinsSerrin,Serrin}).
When $H\equiv 0$  and $\partial\Omega$  is smooth, this curvature condition is that $\partial\Omega$  must have nonnegative mean 
curvature (with respect to the interior normal direction of $\Omega$) at each point (\cite{JenkinsSerrin}). 
However, Leon Simon (\cite{Simon}) has shown that if $\Gamma_{0}\subset \partial\Omega$  is smooth (i.e. $C^{4}$), the mean curvature $\Lambda$  
of $\partial\Omega$  is negative on $\Gamma_{0}$  and $\Gamma$  is a compact subset of $\Gamma_{0},$
then the minimal hypersurface $z=f(x),$  $x\in\Omega,$  extends to $\Omega\cup \Gamma$  as a continuous function, even though $f$  may not equal $\phi$  on $\Gamma.$  
Since \cite{Simon} appeared, the requirement that $H\equiv 0$  has been eliminated and the conclusion remains similar to that which Simon reached 
(see, for example, \cite{Bour, LauLin, Lin}).

How important is the role of boundary smoothness in the conclusions reached in \cite{Simon}?  
We shall show, by constructing suitable domains $\Omega$  and Dirichlet data $\phi,$   
that the existence of a ``nonconvex corner'' $P$  in $\Gamma$  can cause the unique generalized (e.g.  variational) solution  
to be discontinuous at $P$  even if $\Gamma\setminus \{P\}$  is smooth and the generalized mean curvature $\Lambda^{*}$  
(i.e. \cite{Serrin}) of  $\Gamma$  at $P$  is $-\infty$;  this shows that some degree of smoothness of 
$\Gamma$  is required to obtain the conclusions in \cite{Simon}.
We shall prove the following 
\begin{thm}
\label{Main Theorem}
Let $n\in \Natural,$  $n\ge 2,$  and assume there exists $\lambda>0$  such that $|H(x,t)|\le \lambda$  for $x\in\Real^{n}$  and $t\in\Real.$  
Then there exist a domain $\Omega\subset\Real^{n}$  and a point $P\in\partial\Omega$   such that   
\begin{itemize}
\item[(i)] $\partial\Omega\setminus\{P\}$  is smooth ($C^{\infty}$),
\item[(ii)] there is a neighborhood ${\cal N}$  of $P$  such that $\Lambda(x)<0$  for $x\in {\cal N}\cap \partial\Omega\setminus\{P\},$  where $\Lambda$  
is the mean curvature of $\partial\Omega,$  and   
\item[(iii)] $\Lambda^{*}(P)=-\infty,$  where $\Lambda^{*}$  is the generalized mean curvature of $\partial\Omega,$  
\end{itemize}
and there exists Dirichlet boundary data $\phi\in C^{\infty}\left(\Real^{n}\right)$  such that the minimizer $f\in BV(\Omega)$  of 
\begin{equation}
\label{The_Functional} 
J(u)=\int_{\Omega} |Du| + \int_{\Omega} \int_{0}^{u} H(x,t) dt \ dx +\int_{\partial\Omega} |u-\phi| d {\cal H}^{n-1}, \ \ \ \ u\in BV(\Omega),
\end{equation}
exists and satisfies (\ref{ONE-A}), $f\in C^{2}(\Omega)\cap C^{0}\left(\overline{\Omega}\setminus \{P\}\right) \cap L^{\infty}(\Omega),$   
$f\notin C^{0}\left(\overline{\Omega}\right)$  and $f \neq \phi$  in a neighborhood of $P$  in $\partial\Omega.$      
\end{thm}
\vspace{3mm}

Since there are certainly many examples of Dirichlet problems which have continuous solutions even though their  domains fail to satisfy appropriate smoothness 
or boundary curvature conditions (e.g. by restricting to a smaller domain a classical solution of a Dirichlet problem on a larger domain), the question of necessary or sufficient conditions 
for the continuity at $P$  of a generalized solution of a particular Dirichlet problem is of interest and the examples here suggest (to us) 
that a ``Concus-Finn'' type condition might yield necessary conditions for the continuity at $P$  of solutions; see \S \ref{CFcondition}.  

We view this note as being analogous to other articles (e.g. \cite{FinnShi, HuffMcCuan06, HuffMcCuan09, Korevaar}) which enhance our knowledge of the behavior of solutions of boundary 
value problems for prescribed mean curvature equations by constructing and analyzing specific examples.  
One might also compare Theorem \ref{Main Theorem} with the behavior of generalized solutions of (\ref{ONE-A})-(\ref{ONE-B}) when $\partial\Omega\setminus\{P\}$  
is smooth and $|H(x,\phi(x))|\le (n-1)\Lambda(x)$  for  $x\in \partial\Omega\setminus\{P\}$  (e.g. \cite{EL1986, Lan1985, Lan1988}) and with capillary surfaces (e.g. \cite{LS1}).

\section{Nonparametric Minimal Surfaces in $\Real^{3}$}
\label{BLUE}

In this section, we will assume $n=2$  and $H\equiv 0;$  this allows us to use explicit comparison functions and illustrate our general procedure. 
Let $\Omega$  be a bounded, open set in $\Real^{2}$  with locally Lipschitz boundary $\partial\Omega$  such that a point $P$  
lies on $\partial\Omega$  and there exist distinct rays $l^{\pm}$  starting at $P$  such that $\partial\Omega$  is tangent to 
$l^{+}\cup l^{-}$  at $P.$ 
By rotating and translating the domain, we may assume $P=(0,1)$  and there exists a $\sigma\in \left(-\frac{\pi}{2},\frac{\pi}{2}\right)$  
such that 
$l^{-}=\{\left(r\cos(\sigma),1+r\sin(\sigma)\right) : r\ge 0\},$ 
$l^{+}=\{\left(r\cos(\pi-\sigma),1+r\sin(\pi-\sigma)\right) : r\ge 0\}$  and 
\begin{equation}
\label{PIZZA}
\Omega\cap B\left(P,\delta\right) = \{\left(r\cos(\theta),1+r\sin(\theta)\right) : 0<r<\delta, \theta^{-}(r)<\theta<\theta^{+}(r)\}
\end{equation}
for some  $\delta>0$  and functions $\theta^{\pm}\in C^{0}(\left[0,\delta)\right)$  which satisfy $\theta^{-}<\theta^{+},$
$\theta^{-}(0)=\sigma$  and $\theta^{+}(0)=\pi-\sigma$;  here $B\left(P,\delta\right)$  is the open ball in $\Real^{2}$  centered at $P$ 
of radius $\delta.$  If we set $\alpha=\frac{\pi}{2}-\sigma,$  then  $\alpha\in (0,\pi)$  and the angle at $P$  in $\Omega$  of  
$\partial\Omega$  has size $2\alpha.$   
As $\sigma<0$  goes to zero, $2\alpha>\pi$  goes to $\pi$  and the (upper) region   
between $l^{-}$  and $l^{+}$  becomes ``less nonconvex'' and approaches a half-plane   through $P.$
We will show that for each choice of $\sigma\in \left(-\frac{\pi}{2},0\right),$  there is a domain $\Omega$  as above and a choice of 
Dirichlet data $\phi\in C^{\infty}\left(\partial\Omega\right)$  such that the solution of (\ref{ONE-A})-(\ref{ONE-B})  for $\Omega$ 
and $\phi$  is discontinuous at $P.$  

Fix $\sigma\in \left(-\frac{\pi}{2},-\frac{\pi}{4}\right).$  Let $\epsilon$  be a small, fixed parameter, say $\epsilon\in (0,0.5),$  and let 
$a=a(\sigma)\in (1,2)$  be a parameter to be determined.
Set $\tau=(1+\epsilon)\cot(-\sigma)$  and  $r_{1}=\sqrt{\tau^{2}+(1+\epsilon)^{2}}.$  
Define $h_{2/\pi}\in C^{2}((0,2)\times (-1,1))$  by 
\[
h_{2/\pi}(x_1,x_2) =
\frac{2}{\pi}\ln\left(\frac{\cos\left(\frac{\pi x_2}{2}\right)}{\sin\left(\frac{\pi x_1}{2}\right)}
\right).
\]
Notice that the graph of $h_{2/\pi}$  is part of Scherk's first surface, so $\mathrm {div}(Th_{2/\pi})=0$  on $(0,2)\times (-1,1),$ 
and $h_{2/\pi}(t,t-1)=0$  for each $t\in (0,2).$   
A computation using L'Hospital's Rule shows 
\begin{equation}
\label{ZipZoom}
\lim_{t\to 0^{+}} h_{2/\pi}((t\cos(\theta),1+t\sin(\theta)))  =  \frac{2}{\pi} \ln(-\tan(\theta)), \ \ \ \theta\in \left(-\frac{\pi}{2},0\right)
\end{equation}

\begin{figure}
\centering
\includegraphics{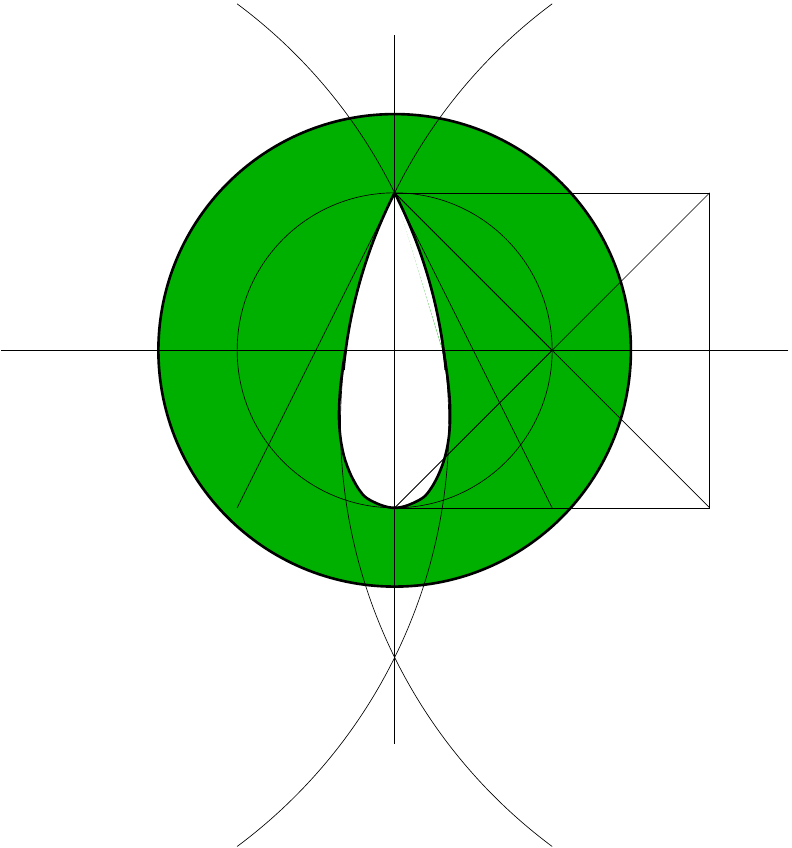}  
\caption{$\Omega$}
\label{FigureOne}
\end{figure}

Let $D=B\left({\cal O},1\right)\cap B\left((\tau,-\epsilon),r_{1}\right)\cap B\left((-\tau,-\epsilon),r_{1}\right)$  be the intersection 
of three open disks  and let  $E\subset D$  be a strictly convex domain such that $\{x\in \partial E:x_{2}< 1\}$  is a $C^{\infty}$  curve, 
$E\cap\{x_{2}\ge 0\} = D\cap \{x_{2}\ge 0\},$  $E$  is symmetric with respect to the $x_{2}-$axis  and  $(0,-1)\in \partial E;$  
here ${\cal O}$  denotes $(0,0).$
Define 
\[
\Omega = B\left({\cal O}, a\right)\setminus \overline{E}
\]
(see Figure \ref{FigureOne}); notice that $P\in \partial\Omega$  and (\ref{PIZZA}) holds with the choice of $\sigma$  above.   
If we set $C=\{(x_{1},x_{2})\in \Real^{2} : 0<x_{1}<1, x_{1}-1<x_{2}<1-x_{1} \},$  then (\ref{ZipZoom}) implies 
$\sup_{x\in C\cap \partial E} h_{2/\pi}(x)<\infty.$

\begin{figure}
\centering
\includegraphics{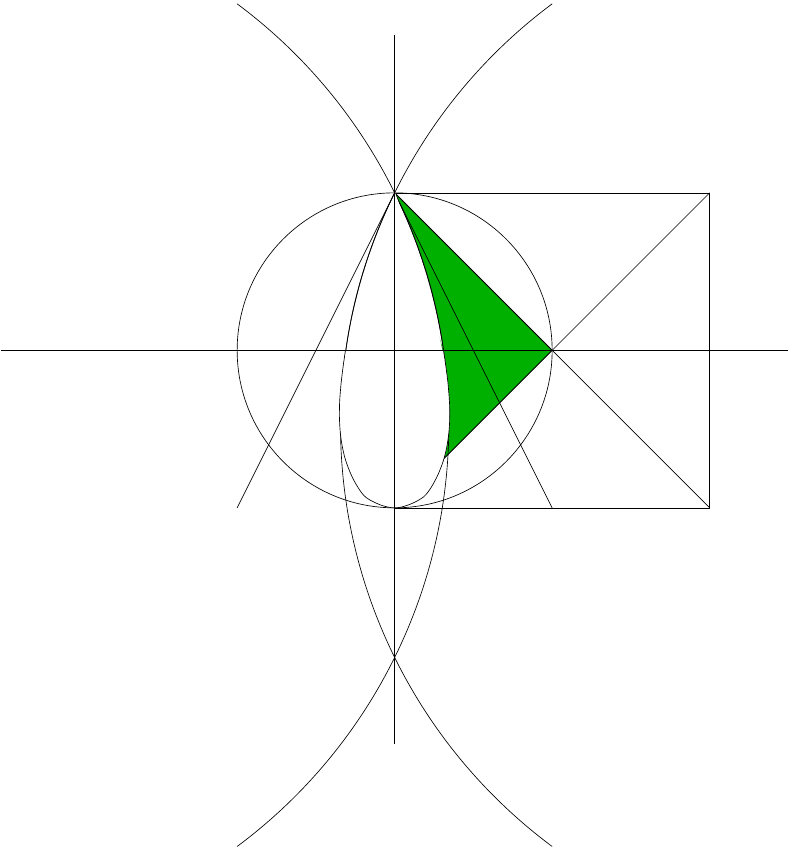}  
\caption{$\Omega\cap C,$  the domain of the comparison function for (\ref{TRex})}
\label{FigureTwo}
\end{figure}

Let $m>\max\{r_{1}\cosh^{-1}\left(\frac{2+\sqrt{\tau^{2}+\epsilon^{2}}}{r_{1}}\right),\sup_{x\in C\cap \partial E} h_{2/\pi}(x)  \}.$  
Notice that $m$  is independent of the parameter $a.$   
Define $\phi\in C^{\infty}\left(\partial\Omega\right)$  by $\phi=0$  on $\partial B\left({\cal O},a\right)$  and $\phi=m$ on $\partial E.$  
Let $f$  be the variational solution of (\ref{ONE-A})-(\ref{ONE-B}) with $\phi$  as given here (e.g. \cite{Ger2,Giu178}).  
Since $\phi\ge 0$  on $\partial\Omega$  and $\phi>0$  on $\partial E,$  $f\ge 0$  in $\Omega$  (e.g. Lemma \ref{Four} (with $h\equiv 0$)) and 
so $f>0$   in $\Omega$  (e.g. the Hopf boundary point lemma). 
Notice that $h_{2/\pi}=0<f$  on $\Omega\cap\partial C$  and $h_{2/\pi}<\phi$  on $C\cap \partial E = C\cap \partial\Omega$  and therefore 
$h_{2/\pi}<f$  on $\Omega\cap C$  (see Figure \ref{FigureTwo}).    
Together with (\ref{ZipZoom}), this implies 
\begin{equation}
\label{TRex}
\liminf_{\Omega\cap C\ni x\to P} f(x) \ge \frac{2}{\pi} \ln(\tan(-\sigma))>0.
\end{equation}

\begin{figure}
\centering
\includegraphics{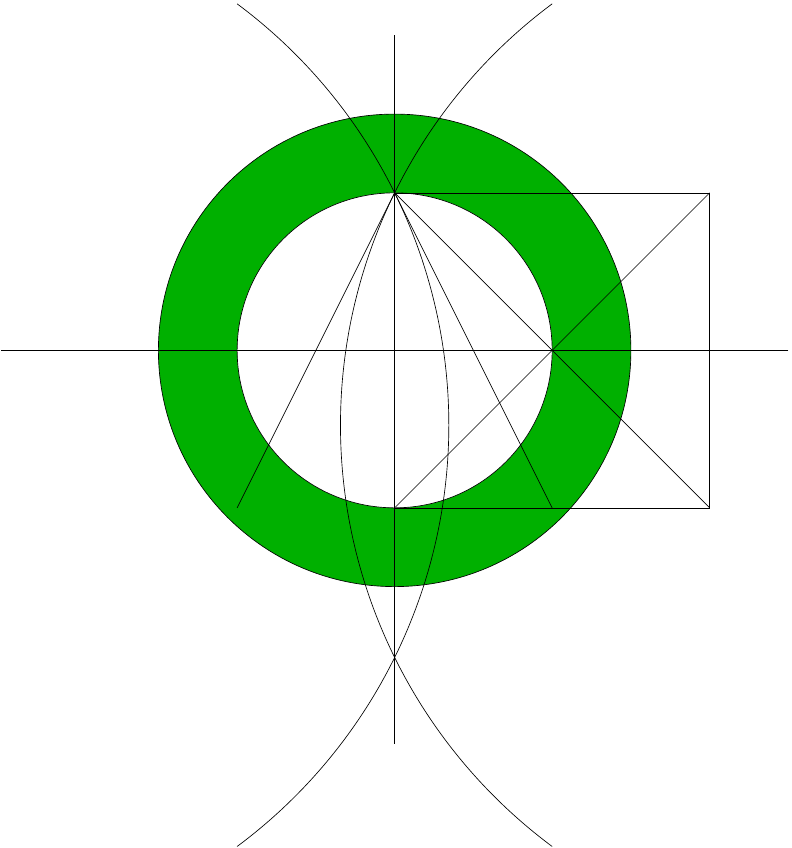}  
\caption{$W,$  the domain  of the comparison function for (\ref{Raptor})}
\label{FigureThree}
\end{figure}

Set  $W=B\left({\cal O},a\right) \setminus \overline{B\left({\cal O},1\right)}$  (see Figure \ref{FigureThree});  then $W\subset \Omega$.  
Define  $g\in C^{\infty}(W)\cap C^{0}(\overline{W})$  by $g(x) = \cosh^{-1}\left(a\right)-\cosh^{-1}\left(|x|\right)$
and notice that the graph of $g$  is part of a catenoid, where $g=0$  on $\partial B\left({\cal O},a\right)$  and 
$g=\cosh^{-1}\left(a\right)$  on $\partial B\left({\cal O},1\right).$  
It follows from the General Comparison Principle (e.g. \cite{FinnBook}, Theorem 5.1)  that $f\le g$  on $W$  and therefore 
\begin{equation}
\label{Raptor}
f\le \cosh^{-1}\left(a\right) \ \ \ \ \  {\rm on} \ \ W.  
\end{equation}
If we select $a>1$  so that $\cosh^{-1}\left(a\right) < \frac{2}{\pi} \ln(\tan(-\sigma)),$  then (\ref{TRex}) and (\ref{Raptor}) imply 
that $f$  cannot be continuous at $P.$  
Notice that \cite{Simon} implies $f\in C^{0}\left(\overline{\Omega}\setminus \{P\}\right).$  

This example illustrates the procedure we shall use in \S \ref{Happy}; a somewhat similar approach was used in \cite{FinnShi,Korevaar,LS1,Serrin}. 
The case when $\sigma\in \left[-\frac{\pi}{4},0\right)$  has a similar proof with the changes that $D$  is the intersection of the open disk 
$B\left({\cal O},1\right)$  with the interiors of two ellipses and a Scherk surface with rhomboidal domain (\cite{NitscheBook}, pp. 70-71) is used 
as a comparison surface to obtain the analog of (\ref{TRex}); the details can be found in \cite{MelinThesis}.

\section{Lemmata}

\begin{lem}
\label{Three}
Let $\Omega$  be a bounded open set in $\Real^{n},$  $n\ge 2,$  with locally Lipschitz boundary  and let $\Gamma$  be an open,  $C^{2}$  
subset of $\partial\Omega.$  Let $\phi \in L^{\infty}(\partial\Omega)\cap C^{1,\beta}(\Gamma).$
Suppose $g\in C^{2}(\Omega)\cap L^{\infty}(\Omega)$ is the variational solution of (\ref{ONE-A})-(\ref{ONE-B}) and $g<\phi$  on $\Gamma.$  
Then $\nu\equiv \frac{(\nabla g,-1)}{\sqrt{1+|\nabla g|^{2}}}\in C^{0}\left(\Omega\cup\Gamma \right)$  and 
$\nu\cdot\eta=1$  on $\Gamma,$
where $\eta(x)\in S^{n-1}$  is the exterior unit normal to $\Gamma$  at $x.$ 
\end{lem}
\vspace{3mm}

\noindent {\bf Proof:}  Since $g$  minimizes the functional $J$  in (\ref{The_Functional})  over $BV\left(\Omega\right),$   $g$  also minimizes the functional 
$K(u)=J(u)-\int_{\Gamma} \phi\ d {\cal H}^{n-1}.$  Notice   
\[
K(u)=\int_{\Omega} |Du| + \int_{\Omega} \int_{0}^{u} H(x,t) dt \ dx +\int_{\partial\Omega\setminus\Gamma} |u-\phi| d {\cal H}^{n-1} - \int_{\Gamma} u\ d {\cal H}^{n-1} 
\]
for each  $u\in BV\left(\Omega\right)$  with $tr(u)\le \phi$  on $\Gamma;$  in particular, this holds when $u=g.$  
Therefore, for each $x\in\Gamma,$  there exists $\rho>0$  such that $\partial\Omega\cap B_{n}(x,\rho)\subset\Gamma$   and 
the Lemma follows as in \cite{KorevaarSimon}.  \qed
\vspace{3mm}

\begin{lem}
\label{Four}
Let $\Omega$  be a bounded open set in $\Real^{n},$  $n\ge 2,$  with locally Lipschitz boundary, $\phi,\psi \in L^{\infty}(\partial\Omega)$  
with $\psi\le\phi$  on $\partial\Omega,$  $H_{0}\in C^{2}\left(\Omega\times\Real\right)$  with $H_{0}(x,t)$  nondecreasing in $t$  for $x\in\Omega,$  
and $H_{0}\ge H$  on $\Omega\times\Real.$    
Consider the boundary value problem 
\begin{eqnarray}
\label{ONE-D}
{\rm div}\left( Tf\right) &  = &  H_{0}(x,f) \ \ \ \ \ {\rm in} \ \ \Omega, \\ 
f & = & \psi \ \ \ \ \ {\rm  on} \ \ \partial\Omega. 
\label{ONE-E}
\end{eqnarray}
Suppose $g\in C^{2}(\Omega)\cap L^{\infty}(\Omega)$ is the variational solution of (\ref{ONE-A})-(\ref{ONE-B}) and either 
(i) $h\in C^{2}(\Omega)\cap L^{\infty}(\Omega)$  is the variational solution of (\ref{ONE-D})-(\ref{ONE-E}) or
(ii) $\psi\in C^{0}(\partial{\Omega}),$  $h\in C^{2}(\Omega)\cap C^{0}(\overline{\Omega})$  and $h$  satisfies (\ref{ONE-D})-(\ref{ONE-E}). 
Then $h\le g$  in $\Omega.$
\end{lem}
\vspace{3mm}

\noindent {\bf Proof:}  Let $A=\{x\in\Omega : h(x)>g(x)\}.$  In case (i), let $f=hI_{\Omega\setminus A}+gI_{A},$  where $I_{B}$  is the characteristic 
function of a set $B;$  then a simple calculation using $J(g)\le J(f)$  shows that $J_{1}(f)\le J_{1}(h)$  and therefore $f=h$  and $A=\emptyset,$ 
where $J_{1}(u)=\int_{\Omega} |Du| + \int_{\Omega} \int_{0}^{u} H_{0}(x,t) dt \ dx  +\int_{\partial\Omega} |u-\psi| d {\cal H}^{n-1},$  $u\in BV(\Omega),$  is the 
functional $h$  minimizes.  Case (ii) follows from Lemma 1 of \cite{Williams}.  \qed

\begin{lem}
\label{Mango}
Let $\Omega \subset \{x\in\Real^{2}\ : \ x_{2}>0\}$ be a bounded open set, $n\in \Natural$  with $n\ge 2$   and  $g\in C^{2}\left(\Omega\right).$  
Set $\tilde \Omega = \{(x_{1},x_{2}\omega)\in \Real^{n}\ : \ (x_{1},x_{2})\in \Omega, \omega\in S^{n-2} \}$  and 
define $\tilde g\in C^{2}\left( \tilde \Omega \right)$  by 
$\tilde g(x_{1},x_{2}\omega)=g(x_{1},x_{2})$  for  $(x_{1},x_{2})\in \Omega, \omega\in S^{n-2}.$
Then, for $x=(x_{1},\dots,x_{n})=(x_{1},r\omega) \in \tilde \Omega$  with $r=\sqrt{x_{2}^{2}+\dots+x_{n}^{2}},$  $\omega=\frac{1}{r}(x_{2},\dots,x_{n})$    
and $(x_{1},r)\in \Omega,$   we have 
\[
{\rm div}\left( \frac{\nabla \tilde g}{\sqrt{1+|\nabla \tilde g|^{2}}}\right)(x) 
= {\rm div}\left( \frac{\nabla g}{\sqrt{1+|\nabla g|^{2}}}\right)(x_{1},r)
+ \frac{n-2}{r}\frac{g_{x_{2}}(x_{1},r)}{\sqrt{1+|\nabla g(x_{1},r)|^{2}}}.
\]
In particular, if $H\ge 0,$  $R>0,$  $\Omega\subset  \{x\in\Real^{2}\ : \ x_{2}\ge R\}$  and 
\[
{\rm div}\left( \frac{\nabla g}{\sqrt{1+|\nabla g|^{2}}}\right)\ge H +     \frac{n-2}{R} \ \ \ \ \  {\rm on} \ \ \Omega,
\]
then ${\rm div}\left( \frac{\nabla \tilde g}{\sqrt{1+|\nabla \tilde g|^{2}}}\right)\ge H$  on $\tilde\Omega$
\end{lem}
\vspace{3mm}

\noindent {\bf Proof:} Notice that $1+|\nabla\tilde g|^{2}=1+|\nabla g|^{2},$
\[
\left( 1+|\nabla \tilde g|^{2}\right) \triangle \tilde g = 
\left( 1+|\nabla g|^{2}\right) \left( \triangle  g  + \frac{n-2}{r} g_{x_{2}}  \right),
\]
\[
\sum_{i,j=1}^{n} \frac{\partial \tilde g}{\partial x_{i}} \frac{\partial \tilde g}{\partial x_{j}}  \frac{\partial^{2} \tilde g}{\partial x_{i}\partial x_{j}} 
= \left( \frac{\partial g}{\partial x_{1}}\right)^{2}  \frac{\partial^{2} g}{\partial x_{1}^{2}} 
+2\frac{\partial g}{\partial x_{1}}\frac{\partial g}{\partial x_{2}}\frac{\partial^{2} g}{\partial x_{1}\partial x_{2}} 
+ \left( \frac{\partial g}{\partial x_{2}}\right)^{2}  \frac{\partial^{2} g}{\partial x_{2}^{2}} 
\] 
and so 
\[
\left( 1+|\nabla \tilde g|^{2}\right) \triangle \tilde g - 
\sum_{i,j=1}^{n} \frac{\partial \tilde g}{\partial x_{i}} \frac{\partial \tilde g}{\partial x_{j}}  \frac{\partial^{2} \tilde g}{\partial x_{i}\partial x_{j}} 
\]
\[
 =  \left(1+g_{x_{2}}^{2}\right)g_{x_{1}x_{1}} - 2g_{x_{1}}g_{x_{2}}g_{x_{1}x_{2}}+\left(1+g_{x_{1}}^{2}\right)g_{x_{2}x_{2}} 
 + \frac{n-2}{r} \left(1+g_{x_{1}}^{2}+g_{x_{2}}^{2} \right) g_{x_{2}}.
\]
The lemma follows from this.  \qed

\section{The $n-$dimensional case}
\label{Happy}

Let $B_{k}\left(x,r\right)$  denote the open ball in $\Real^{k}$  centered at $x\in \Real^{k}$  with radius $r>0$  and ${\cal O}_{k}=(0,\dots,0)\in \Real^{k},$  for $k\in\Natural.$  
Now consider $n\ge 2$  and set  
\[
\lambda=\sup_{(x,t)\in \Real^{n}\times\Real}|H(x,t)|;
\]  
if $\lambda=0,$  replace it with a positive constant.  For each $a\in \left(0,\frac{n}{\lambda}\right)$  and $Q\in \Real^{n},$  we have 
\begin{equation}
\label{BOOK}
\int_{B_{n}(Q,a)} \lambda^{n} dx < n^{n}\omega_{n}.   
\end{equation}
By translating our problem in $\Real^{n},$  we may (and will) assume $Q={\cal O}_{n}.$
By Proposition 1.1 and Theorem 2.1 of \cite{Giu176}, we see that if $\Omega$  is a bounded, connected, open set in  $\Real^{n}$  with Lipschitz-continuous boundary,  
$\overline{\Omega}\subset B_{n}\left({\cal O}_{n},\frac{n}{\lambda}\right)$  and $\phi\in L^{1}\left(\partial\Omega\right),$   then the functional $J$  in (\ref{The_Functional}) has a minimizer 
$f\in BV(\Omega),$   $f\in C^{2}(\Omega)$  and $f$  satisfies (\ref{ONE-A}).  

The proof in \S \ref{CDO} consists of setting some parameters (e.g. $p,$  $r_{1},$  $r_{2},$  $m_{0},$  $b,$  $c,$   $\tau,$  $\sigma,$  $a$), determining the domain $\Omega,$  
finding different comparison functions (e.g. $g_{1},$  $g^{[u]},$  $k_{\pm},$  $k_{2},$  $k_{3},$  $k_{4}$), and mimicking (\ref{TRex}) and (\ref{Raptor}) 
to show that the variational solution $f$  of (\ref{ONE-A})-(\ref{ONE-B}) is discontinuous at a nonconvex corner.  
In particular, we use a torus (i.e. $j_{a}$) to obtain (\ref{Hopping}), unduloids (i.e. $k_{\pm},$  $k_{2}$)  to obtain (\ref{CatX}) (an analog of (\ref{Raptor})) and 
nodoids (i.e. $g_{1},$  $g^{[u]}$), unduloids (i.e. $k_{\pm},$  $k_{4}$)  and a helicoidal function (i.e. $h_{2}$)  to obtain (\ref{Cab})  (an analog of (\ref{TRex}))   
and prove that $f$  is discontinuous at $P=(0,p,0,\dots,0)\in \Real^{n}\in \partial\Omega.$

\subsection{Codimension $1$  singular set}
\label{CDO}

In this section, we will obtain a domain $\Omega$  as above and $\phi\in C^{\infty}(\Real^{n})$  such that $P\in\partial\Omega,$   the minimizer $f$  of 
(\ref{The_Functional}) is discontinuous at $P,$   $\partial\Omega\setminus T$  is smooth ($C^{\infty}$)
and $f\in C^{2}(\Omega)\cap C^{0}\left(\overline{\Omega}\setminus T\right),$  where $T$  is a smooth set of dimension $n-2$  (i.e. $T$  has codimension $1$  in $\partial\Omega$).  
We will use portions of nodoids,  unduloids and  helicoidal surfaces with constant mean curvature as comparison functions.  
For the convenience of the reader, we will denote functions whose graphs are subsets of nodoids with the letter $g$  
(e.g. $g_{1}(x_{1},x_{2})$), subsets of CMC helicoids with the letter $h$  and subsets of unduloids (or onduloids) with the letter $k.$   

Let ${\cal N}_{1}\subset \Real^{3}$  be a nodoid which is symmetric with respect the $x_{3}$-axis  and has mean curvature $1$  (when ${\cal N}_{1}$  is oriented ``inward'', 
so that the unit normal $\vec N_{{\cal N}_{1}}$  to ${\cal N}_{1}$  points toward the $x_{3}$-axis at the points of ${\cal N}_{1}$  which are furthest from the $x_{3}$-axis).  
Let $s_{1}=\inf_{(x,t)\in {\cal N}_{1}} |x|$  be the inner neck size of ${\cal N}_{1}$  and let $s_{3}$  satisfy the condition that the unit normal to ${\cal N}_{1}$  is vertical 
(i.e. parallel to the $x_{3}$-axis) at each point $(x,t)\in \Real^{2}\times\Real$  of ${\cal N}_{1}$  at which $|x|=s_{3};$  then $s_{1}<s_{3}.$   Let $s_{2}\in (s_{1},s_{3}).$    
(Notice that we can assume $s_{2}/s_{1}$  is close to $s_{3}/s_{1}$  if we wish.)  

Let us fix  $0<p<\frac{1}{\lambda}$  and set $w=(0,p)\in\Real^{2},$   $P=(0,p,0,\dots,0)\in \Real^{n}.$  
Let $m_{0}=\lambda/2+(n-2)/(p/3).$   We shall assume $r_{2}=s_{2}/m_{0}<p/3;$  if necessary, we increase $m_{0}$  to accomplish this.  
Let $r_{1}=s_{1}/m_{0}$  and $r_{3}=s_{3}/m_{0}.$  
Let ${\cal N}=\{(m_{0})^{-1}X\in \Real^{3} : X\in {\cal N}_{1} \};$  then ${\cal N}$  is a nodoid with mean curvature $m_{0}.$  
Set $\Delta_{1}=\{x\in\Real^{2} : r_{1}<|x|< r_{2} \}.$  
Fix $b\in \left(0,\frac{1}{4m_{0}}\left(1+2m_{0}p-\sqrt{1+4m_{0}^{2}p^{2}}\right)\right).$

Define $g_{1}\in C^{\infty}\left(\Delta_{1}\right)\cap C^{0}\left(\overline{\Delta_{1}}\right)$  to be a function whose graph is a 
subset of ${\cal N}$  on which  $\vec N_{\cal N}=(n_{1},n_{2},n_{3})$  satisfies $n_{3}\ge 0;$  
then 
\begin{equation}
\label{STAR}
{\rm div}\left( \frac{\nabla g_{1}}{\sqrt{1+|\nabla g_{1}|^{2}}}\right) = m_{0}\ge \lambda+\frac{2(n-2)}{p/3}. 
\end{equation}
By moving ${\cal N}$  vertically, we may assume $g_{1}(x)=0$  when $|x|=r_{2};$  then $g_{1}>0$  in $\Delta_{1}.$   
Notice that $\frac{\partial g_{1}}{\partial x_{1}}(r_{1},0)=-\infty$  and $\frac{\partial g_{1}}{\partial x_{1}}(r_{2},0)<0$;
then there exists a $\beta_{0}>0$  such that,  for each $\theta\in \Real,$  
\begin{equation}
\label{Dog}
\frac{\partial}{\partial r} \left( g_{1}(r\Theta) \right)<-\beta_{0}  \ \ \ \ \ {\rm for} \ \ r_{1}<r<r_{2}, 
\end{equation}
where $\Theta=(\cos(\theta),\sin(\theta)).$   Fix $\beta\in (0,\beta_{0}).$  
Let 
\begin{equation}
\label{Shockers}
0<\tau <  \min \left\{\frac{pr_{1}}{\sqrt{r_{2}^{2}-r_{1}^{2}}},  \frac{2(1-p\lambda)}{\lambda(2-p\lambda)}, \frac{b(4p-b)}{4(2p-b)} \right\}.
\end{equation}

Consider $\sigma\in \left(-\frac{\pi}{2},0\right).$  
Notice that the distance between  $L$  and the point $(0,p-r_{2})$  is $r_{2}\cos(\sigma),$  where $L$  is the closed sector given by 
\[
L=\{ \left(r\cos(\theta),p+r\sin(\theta)\right) : r\ge 0, \sigma\le\theta\le\pi-\sigma \}.
\]  
Define $r_{4}=\sqrt{p^{2}+\tau^{2}}$  and   
\[
M= B_{2}\left((\tau,0),r_{4}\right) \cap B_{2}\left((-\tau,0),r_{4}\right).
\]
Notice that $\tau<\frac{b(4p-b)}{4(2p-b)}$  and therefore $B_{2}\left({\cal O}_{2},\frac{a+p}{2}-b \right)\subset M$  if $p<a<p+b.$  

\begin{figure}
\centering
\includegraphics{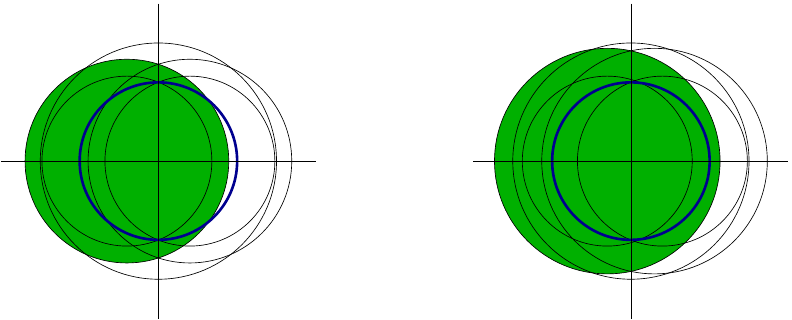}  
\caption{(a) $B_{2}\left({\cal O}_{2},p\right)\nsubseteq B_{2}\left((-\tau,0),\frac{2}{\lambda}-r_{4}\right)$  
(b) $B_{2}\left({\cal O}_{2},p\right)\subset B_{2}\left((-\tau,0),\frac{2}{\lambda}-r_{4}\right)$ }
\label{FigureEight}
\end{figure}

Set $\sigma=-\arctan(\tau/p);$  then $\cos(\sigma)>\frac{r_{1}}{r_{2}},$  since $\tau<\frac{p\sqrt{r_{2}^{2}-r_{1}^{2}}}{r_{1}},$  
and  $L\cap \overline{B_{2}}=\emptyset,$  where $B_{2}=B_{2} \left( (0,p-r_{2}) ,r_{1} \right).$   Therefore there exists a $\delta_{1}>0$  such that if 
$u=(u_{1},u_{2})\in \partial B_{2}({\cal O}_{2},p)$   with $|u-w|<\delta_{1},$  then  
\begin{equation}
\label{Cat2}
B_{2}\left(\frac{p-r_{2}}{p}u,r_{1}\right)  \subset M.  
\end{equation}
Since $\tau<\frac{2(1-p\lambda)}{\lambda(2-p\lambda)},$  we have $\tau-\left( \frac{2}{\lambda}-r_{4}\right)< -p$    and so 
$B_{2}\left({\cal O}_{2},p\right)\subset B_{2}\left((\tau,0),\frac{2}{\lambda}-r_{4}\right)$  (see Figure \ref{FigureEight} (b)).  
Notice that  
\begin{equation}
\label{SALAD}
M\setminus \{(0,\pm p)\} = \{\left(r\cos(\theta),p+r\sin(\theta)\right) : 0<r<2p, \theta^{-}(r)<\theta<\theta^{+}(r)\}
\end{equation}
for some  functions $\theta^{\pm}\in C^{0}(\left[0,\delta)\right)$  which satisfy $\theta^{-}<\theta^{+},$
$\theta^{-}(0)=-\pi-\sigma$  and $\theta^{+}(0)=\sigma.$

Let $a>p$  and set 
\[
{\cal T} = \left\{ \left(\left(\frac{a+p}{2}+b\cos v\right)\cos u,\left(\frac{a+p}{2}+b\cos v\right)\sin u,b\sin v+c \right) : (u,v)\in R \right\},
\]
where $R=[0,2\pi]\times [-\pi,0]$  and  $0<c<b;$  since $b< \frac{1}{4m_{0}}\left(1+2m_{0}p-\sqrt{1+4m_{0}^{2}p^{2}}\right),$  we see that 
$\frac{(a+p)/2-2b}{4b((a+p)/2-b)}>m_{0}$  for all $a\ge p.$  
We shall assume 
\begin{equation}
\label{Assumption_a}
a\in \left(p,\min\{p+b,1/\lambda\}\right)
\end{equation}
and $c=\sqrt{b^{2}-\left(\frac{a-p}{2}\right)^{2}}.$  
Notice that ${\cal T}$  is the lower half of a torus whose mean curvature (i.e. one half of the trace of the shape operator) at each point is greater than $m_{0}.$  
Let ${\cal T}$  be the graph of a function $j_{a}$  over $\Delta_{a}=\{x\in \Real^{2} : \frac{a+p}{2}-b\le |x| \le \frac{a+p}{2}+b\};$  then $j_{a}(x)=0$  on $|x|=a$  and $|x|=p,$  $j_{a}(x)<0$  on $p<|x|<a$  and 
$j_{a}(x)>0$  on $\frac{a+p}{2}-b\le |x|<p$  and $a<|x|\le \frac{a+p}{2}+b$  for $x\in \Real^{2}.$    
Notice that $|j_{a}(x)|< \frac{1}{2m_{0}}$  for all $x\in \Delta_{a}.$

\begin{figure}[h]
\centering
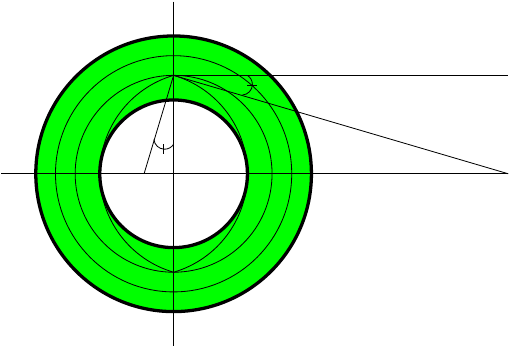
\caption{The domain of $j_{a}$}
\label{FigureFour}
\end{figure}

\begin{figure}
\centering
\includegraphics{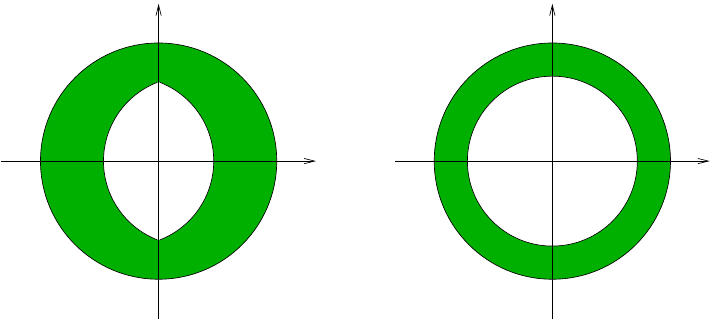}  
\caption{(a)  $\Pi_{1,j}\left(\Omega\right)$  for $2\le j\le n$  (b) $\Pi_{i,j}\left(\Omega\right)$  for  $2\le i<j\le n$}
\label{FigureFive}
\end{figure}

Set 
\begin{equation}
\label{AAAA}
\Omega=B_{n}\left({\cal O}_{n},a\right)\setminus \overline{\cal M},
\end{equation}
where ${\cal M}=\tilde M = \{(x_{1},x_{2}\omega)\in \Real^{n}\ : \ (x_{1},x_{2})\in M, \omega\in S^{n-2} \}.$ 
If we define $\Pi_{i,j}(A)=\{(x_{i},x_{j}):(x_{1},\dots,x_{n})\in A,\ x_{k}=0 \ \ {\rm for}\ \ k\neq i,j\}$  for $A\subset \Real^{n}$  and $1\le i<j\le n,$  then 
$\Pi_{1,j}(\Omega)=B_{2}\left({\cal O}_{2},a\right)\setminus \overline{M}$  for $2\le j\le n$   and 
$\Pi_{i,j}(\Omega)=B_{2}\left({\cal O}_{2},a\right)\setminus \overline{B_{2}\left({\cal O}_{2},1\right)}$  for $2\le i<j\le n$  (see Figure \ref{FigureFive}).  

\begin{figure}
\centering
\includegraphics{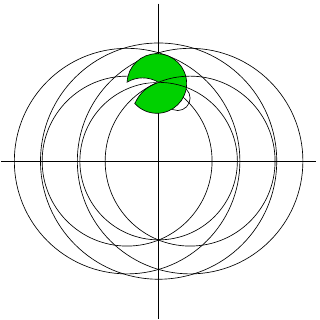}  
\caption{${\cal R}$}
\label{FigureSix}
\end{figure}

We wish to select a helicoidal surface in $\Real^{3}$   (e.g. \cite{DC_MD}) with constant mean curvature $m_{0},$  axis $\{w\}\times\Real$  
and pitch $-\beta$  (recall $-\beta \in (-\beta_{0},0))$, which we will denote ${\cal S};$  then, for each $t\in\Real,$  
$k_{t}\left({\cal S}\right)={\cal S},$   where $k_{t}:\Real^{3}\to \Real^{3}$  is the helicoidal motion  given by 
$
k_{t}(x_{1},x_{2},x_{3})= (l_{t}(x_{1},x_{2}),x_{3}-\beta t)
$
with  $l_{t}:\Real^{2}\to \Real^{2}$   given by 
\[
l_{t}(x_{1},x_{2})=(x_{1}\cos(t)+(x_{2}-p)\sin(t),p-x_{1}\sin(t)+(x_{2}-p)\cos(t)).
\]
Set $c_{0}=\frac{1}{4}\beta \sigma <0;$  by vertically translating ${\cal S},$  we may assume that there is an open  
$c_{0}-$level curve ${\cal L}_{0}$  of ${\cal S}$  with endpoints  $w=(0,p)$  and $b=(b_{1},b_{2})$  such that 
${\cal L}_{0}\subset (0,\infty)\times \Real,$  
${\cal L}=\overline{{\cal L}_{0}}$  is tangent to the (horizontal) line $\Real\times\{p\}$  at $w$  and the slope $m_{v}$  of the tangent line 
to ${\cal L}$  at $v$  satisfies $|m_{v}|<\tan\left(-\sigma/5\right)$  for each $v\in {\cal L}_{0};$
then ${\cal L}\times \{c_{0}\}\subset {\cal S}$  and 
the curves $l_{t}\left({\cal L}_{0}\right),$  $-\frac{7\pi}{8}<t<\frac{7\pi}{8},$  are mutually disjoint.
Notice that the set 
\[
{\cal R}= \{  l_{t}\left({\cal L}_{0}\right) : -\frac{7\pi}{8}<t<\frac{7\pi}{8} \} 
= \bigcup_{-\frac{7\pi}{8}<t<\frac{7\pi}{8}} l_{t}\left({\cal L}_{0}\right)
\]
is an open subset of $\Real^{2}\setminus \left( (-\infty,0]\times \{p\}\right)$  (see Figure \ref{FigureSix}), 
$w\in \overline{{\cal R}}$  and ${\cal S}$  
implicitly defines the smooth function $h_{2}$   on ${\cal R}$  given by $h_{2}(x)= \frac{\beta}{4} (\sigma-4t)$  if 
$x\in l_{t}\left({\cal L}_{0}\right)$  for some $t\in (-\pi/2,\pi/2).$  
Notice that $B_{2}\left(w,b_{1}\right) \cap \{x_{1}>0\} \subset {\cal R}.$  
Now $l_{t}\left({\cal L}_{0}\right)\cap M = \emptyset$  for $t\in \left(3\sigma/4,\sigma/4\right)$  and, by making $b_{1}>0$  
sufficiently small, we may  assume that 
\begin{equation}
\label{STARSHINE}
l_{t}\left({\cal L}_{0}\right)\subset B_{2}({\cal O}_{2},p)\setminus M  \ \ \ \ \ {\rm for \ \ each} \ \ \ \ \  t\in \left(3\sigma/4,\sigma/4\right).
\end{equation}  
Notice that $h_{2}<\frac{\beta(2\sigma^{2}-\pi)}{8\sigma}$  on $l_{t}\left({\cal L}_{0}\right)$  for $-\frac{\pi}{2}<t<\frac{7\pi}{8}.$ 

\begin{figure}
\centering
\includegraphics[trim = 0mm 22mm 0mm 0mm, clip, width=8cm]{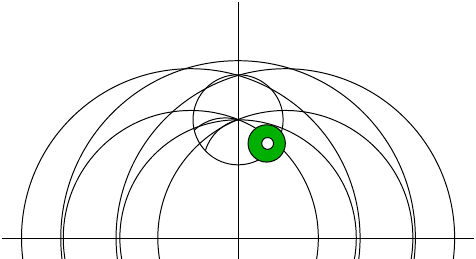} 
\caption{${\cal D}^{[u]};$   $\Omega\cap \tilde{\cal D}^{[u]}$  is the domain  of the comparison function for (\ref{Bingo})}
\label{FigureSeven}
\end{figure}

Let us fix $u=(u_{1},u_{2})\in \partial B_{2}({\cal O}_{2},p)$  such that  $|u-w|<\min\{\delta_{1},b_{1}\}$  and $u_{1}>0.$
Then there exists $\theta_{u}\in (0,\pi/2)$  such that $u=(p\cos(\theta_{u}),p\sin(\theta_{u})).$
Define $g^{[u]}(x)=g_{1}\left(x+\frac{r_{2}-p}{p}u\right)$  and notice that $g^{[u]}(u)=g_{1}\left(\frac{r_{2}}{p}u\right)=0,$  
since $|\frac{r_{2}}{p}u|=r_{2}.$  Note that the domain 
\[
{\cal D}^{[u]}= 
\{x+\frac{p-r_{2}}{p}u : x\in \Delta_{1}\} = B_{2}\left(\frac{p-r_{2}}{p}u,r_{2}\right) \setminus \overline{ B_{2}\left(\frac{p-r_{2}}{p}u,r_{1}\right) }
\]  
of $g^{[u]}$  is contained in $B_{2}({\cal O}_{2},p)$  since $\partial B_{2}\left(\frac{p-r_{2}}{p}u,r_{2}\right)$  and 
$\partial  B_{2}({\cal O}_{2},p)$  are tangent circles at $u$  and $r_{2}<p$  (see Figure \ref{FigureSeven}). 
Notice that 
\begin{equation}
\label{SPACESHIP}
h_{2}(r\cos(\theta_{u}),r\sin(\theta_{u})) < g^{[u]}(r\cos(\theta_{u}),r\sin(\theta_{u}))   
\end{equation}  
when $p-r_{2}+r_{1}\le r\le p,$  because $h_{2}(u) < 0 = g^{[u]}(u),$  $\beta<\beta_{0}$  and (\ref{Dog}) holds. 

Let 
\[
{\cal N}_{\pm}\subset \{x\in\Real^{2} : r_{4}\le |(x_{1}  \pm \tau,x_{2})| \le \frac{2}{\lambda}-r_{4}\}\times\Real
\]  
be  unduloids in $\Real^{3}$  with mean curvature $\lambda/2$  such that $\{(\mp \tau,0)\}\times\Real$   are the respective axes of symmetry; 
the minimum and maximum radii (or ``neck'' and ``waist'' sizes) of both unduloids  are $r_{4}$  and $\frac{2}{\lambda}-r_{4}$  respectively.  
Set 
\[
\Delta_{\pm}=B_{2}\left( (\mp \tau,0),\frac{2}{\lambda}-r_{4}\right) \setminus \overline{B_{2}\left( (\mp \tau,0),r_{4}\right)}
\]  
and define 
$k_{\pm}\in C^{\infty}\left(\Delta_{\pm}\right)$  so that the graphs of $k_{\pm}$  are subsets of ${\cal N}_{\pm}$  respectively,   
\[
{\rm div}\left( Tk_{\pm}\right)=-\lambda \ \ \ \ \ {\rm in} \ \  \Delta_{\pm},
\]
$\frac{\partial}{\partial r} \left(k_{\pm}\left((\mp p,0)+r\Theta\right) \right)|_{r=r_{4}} = -\infty$  and 
$\frac{\partial}{\partial r} \left(k_{\pm}\left((\mp p,0)+r\Theta\right) \right)|_{r=\frac{2}{\lambda}-r_{4}} = -\infty$ 
for each $\theta\in\Real,$  where $\Theta=(\cos(\theta),\sin(\theta)).$  
We may vertically translate ${\cal N}_{\pm}$  so that $k_{\pm}(x)=0$  for $x\in\Real^{2}$  with $|(x_{1}\pm \tau,x_{2})|=\frac{2}{\lambda}-r_{4}.$  
Notice that $k_{+}\left(0,p\right)=k_{-}\left(0,p\right)=\sup_{\Delta_{+}} k_{+}=\sup_{\Delta_{-}} k_{-}.$  

Let ${\cal N}\subset \{x\in\Real^{2} : p\le |x| \le \frac{2}{\lambda}-p\}\times\Real$  be an unduloid with mean curvature $\lambda/2$  
such that the $x_{3}-$axis is the axis of symmetry and the minimum and maximum radii (or ``neck'' and ``waist'' sizes)  are $p$  and 
$\frac{2}{\lambda}-p$  respectively.
Set $\Delta_{2}=B_{2}\left({\cal O}_{2},\frac{2}{\lambda}-p\right) \setminus \overline{B_{2}\left({\cal O}_{2},p\right)}$  and define 
$k_{2}\in C^{\infty}\left(\Delta_{2}\right)$  so that the graph of $k_{2}$  is a subset of ${\cal N},$  
${\rm div}\left( Tk_{2}\right)=-\lambda$  in  $\Delta_{2},$
$\frac{\partial}{\partial r} \left(k_{2}\left(r\Theta\right) \right)|_{r=p} = -\infty$  and 
$\frac{\partial}{\partial r} \left(k_{2}\left(r\Theta\right) \right)|_{r=\frac{2}{\lambda}-p} = -\infty$ 
for each $\theta\in\Real,$  where $\Theta=(\cos(\theta),\sin(\theta)).$  
\begin{figure}
\centering
\includegraphics{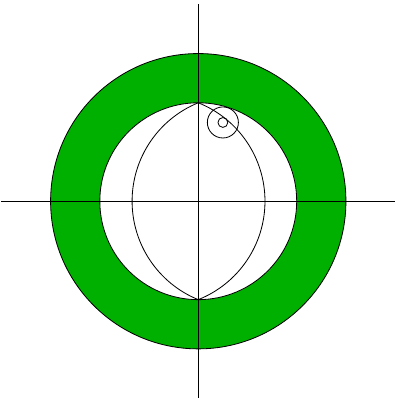}  
\caption{$B_{2}\left({\cal O}_{2},a\right) \setminus \overline{B_{2}\left({\cal O}_{2},p\right)}:$   (\ref{Once upon a time})}
\label{FigureNine}
\end{figure}

Define $\phi\in C^{\infty}\left(\Real^{n}\right)$  so that $\phi=0$  on $\partial B_{n}\left({\cal O}_{n},a\right)$  
and $\phi=m$    on $\partial {\cal M},$  where 
\begin{equation}
\label{Early}
m>\max\{g_{1}(0,r_{1}),\frac{1}{2m_{0}},k_{+}(0,r_{4}-\tau)+ k_{2}(0,p)-k_{2}\left(0,\frac{2}{\lambda}-p\right) \}; 
\end{equation}  
recall then that $m>j_{a}\left(\frac{a+p}{2}-b\right).$  
Let $f$  be the variational solution of (\ref{ONE-A})-(\ref{ONE-B}) with $\Omega$  and $\phi$  as given here;  that is, let $f$  minimize the 
functional given in (\ref{The_Functional}) and notice that the existence of $f$  follows from  (\ref{BOOK}), (\ref{Assumption_a}),  
\S 1.D. of \cite{Giu176}  and \cite{Ger2,Giu178}.
(Notice that there exists $w:B_{2}({\cal O}_{2},a)\setminus M \to \Real$  such that $f= \tilde w.$)  
The comparison principle implies $j_{a}(x)\le f(x)$  for $x\in \Omega$  and so $f(x)\ge j_{a}(x)\ge 0$  if $x\in \Omega$  with $|x|\le p$  
(recall (\ref{Assumption_a}) holds).   In particular, 
\begin{equation}
\label{Hopping}
f(x)\ge 0 \ \ \ \ \ {\rm when} \ \ \  x\in \Omega \ \ {\rm with} \ \  |x| \le p.
\end{equation}

\begin{figure}
\centering
\includegraphics{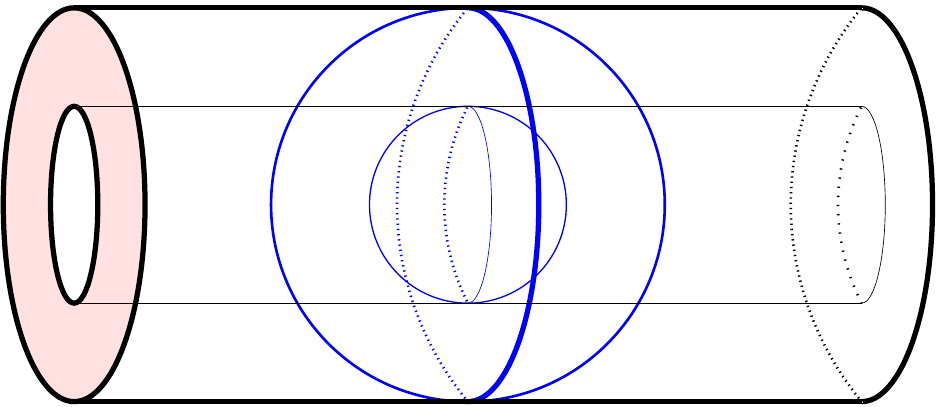}  
\caption{(\ref{Bedtime_Stories}):   $W$  and $B_{n}\left({\cal O}_{n},a\right) \setminus \overline{B_{n}\left({\cal O}_{n},p\right)}$  when $n=3$}
\label{FigureTen}
\end{figure}

Set  $W=\left(B_{2}\left({\cal O}_{2},a\right) \setminus \overline{B_{2}\left({\cal O}_{2},p\right)}\right)\times \Real^{n-2}.$   Now 
\[
\Omega\subset B_{2}\left({\cal O}_{2},a\right) \times \Real^{n-2} \subset B_{2}\left({\cal O}_{2},\frac{2}{\lambda}-p\right) \times \Real^{n-2}
\]
(see Figure \ref{FigureNine}).  Define  $k_{3}(x)=k_{2}\left(x_{1},x_{2}\right)-k_{2}(0,a)$  for $x=(x_{1},x_{2},\dots,x_{n})\in W.$  
Notice that $f=0\le k_{3}$  on $\overline{W}\cap \partial B_{n}\left({\cal O}_{n},a\right),$  
\[
{\rm div}\left( Tf\right)=H(x,f(x))\ge -\lambda= {\rm div}\left( Tk_{3}\right) \ \ \ \ \ {\rm in} \ \ \Omega\cap W
\] 
and $\frac{\partial}{\partial r} \left(k_{2}\left(r\Theta\right) \right)|_{r=p} = -\infty$   
(so that $\lim_{W\ni y\to x} Tk_{3}(y) \cdot \xi(x)=1$  for $x\in \partial B_{2}\left({\cal O}_{2},p\right)\times \Real^{n-2},$   
where $\xi$  is the unit exterior normal to $\partial W$). 
The general comparison principle (e.g. \cite{FinnBook}, Theorem 5.1) then implies 
\begin{equation}
\label{Once upon a time}
f\le k_{3} \ \ \ \ \ {\rm in} \ \ \ \ \  \Omega\cap W
\end{equation}
and, in particular, 
\begin{equation}
\label{Bedtime_Stories}
\limsup_{\Omega\cap W\ni y\to x} f(y) \le k_{3}(x) \ \ \ \ \ {\rm for} \ \ x\in \partial \Omega \cap \overline{W}
\end{equation}
(see Figure \ref{FigureTen}).  
By rotating the axis of symmetry of $W$  through all lines in $\Real^{n}$  containing ${\cal O}_{n}$  (or, equivalently, keeping $W$ 
fixed and rotating $\Omega$  about ${\cal O}_{n}$),   we see that  
\begin{equation}
\label{CatX}
\sup  \{f(x) : x\in B_{n}\left({\cal O}_{n},a\right) \setminus \overline{B_{n}\left({\cal O}_{n},p\right)} \} \le k_{2}(0,p)-k_{2}(0,a).  
\end{equation}  

Now define $k_{4}\in C^{\infty}\left(\Delta_{+}\times \Real^{n-2}\right)\cap  C^{0}\left(\overline{\Delta_{+}}\times \Real^{n-2}\right)$  
by  
\[
k_{4}(x)=k_{+}(x_{1},x_{2}) + k_{2}(0,p)-k_{2}(0,a), \ \ \ x=(x_{1},x_{2},\dots,x_{n})\in \overline{\Delta_{+}}\times \Real^{n-2}.
\]
Combining (\ref{ONE-A}) and (\ref{CatX}) with the facts that ${\rm div}\left(Tk_{4}\right)=-\lambda$   in $\Delta_{+}\times \Real^{n-2}$  and  
$\lim_{\Delta_{+}\times \Real^{n-2}\ni y\to x} Tk_{4}(y) \cdot \xi_{+}(x)=1$  for $x\in \partial B_{2}\left( (-\tau,0),r_{4}\right)\times \Real^{n-2},$  
where $\xi_{+}$  is the inward unit normal to $\partial B_{2}\left( (-\tau,0),r_{4}\right)\times \Real^{n-2},$ 
we see that 
\begin{equation}
\label{CatY}
f\le k_{4} \ \ \ {\rm in} \ \ \Omega \cap \left(\Delta_{+}\times \Real^{n-2}\right).   
\end{equation}
(If Figure \ref{FigureEight} (a) held, then (\ref{CatY}) would not be valid.)  
Now let $L:\Real^{n}\to \Real^{n}$  be any rotation about ${\cal O}_{n}$  which satisfies $L(\Omega)=\Omega,$  
notice that $f\circ L$  satisfies (\ref{ONE-A})-(\ref{ONE-B}) and apply the previous argument 
to obtain $f\circ L \le k_{4}$  in $\Omega \cap \left(\Delta_{+}\times \Real^{n-2}\right)$  and therefore  
\begin{equation}
\label{Taxi}
\sup  \{f(x) : x\in \partial {\cal M}  \}  \le k_{4}(p,0) <m.
\end{equation}
\begin{figure}
\centering
\includegraphics{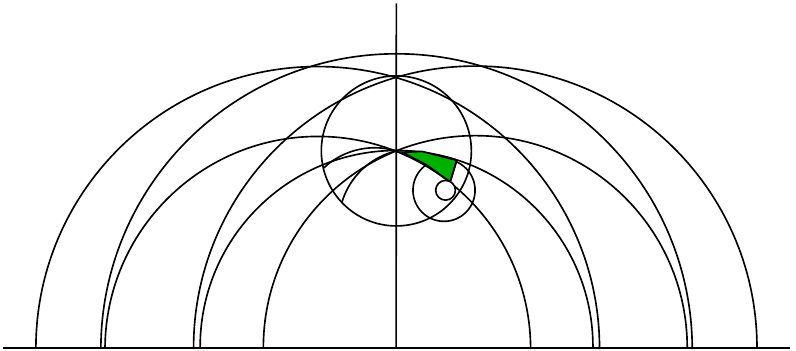}  
\caption{$A:$   (\ref{Jerry})}
\label{FigureEleven}
\end{figure}
From Lemma \ref{Three}, we see that the downward unit normal to the graph of $f,$  $N_{f},$  satisfies
$N_{f} = (\nu,0)$  on $\partial {\cal M} \setminus \{ (0,p\omega) : \omega\in S^{n-2} \}$  and 
\begin{equation}
\label{Doggy}
\lim_{\Omega\ni y\to x} Tf(y)\cdot \nu(x) =1   \ \ \ \ \ {\rm for} \ \ x\in\partial {\cal M} \setminus \{ (0,p\omega) : \omega\in S^{n-2} \}.
\end{equation}
Let us write $B=B_{2}\left(\frac{p-r_{2}}{p}u,r_{2}\right);$   then $\tilde g^{[u]}=0 \le f$  on $\Omega\cap\partial \tilde B$  and 
$\tilde g^{[u]}\le g_{1}(r_{1},0)<\phi$  on $\tilde B \cap \partial M.$   
It follows from (\ref{ONE-A}), (\ref{STAR}) and Lemma \ref{Mango} that 
\begin{equation}
\label{Bingo}
\tilde g^{[u]}<f \ \ \ \ \  {\rm on} \ \ \ \ \ \Omega\cap  \tilde {\cal D}^{[u]}=\Omega \cap \tilde B.
\end{equation}

Set $U=\{r\left(\cos(\theta),\sin(\theta)\omega\right)\in\Omega : r\in (0,p), \theta\in (0,\theta_{u}),\omega\in S^{n-2} \}.$
If we write $\partial_{1}U=\{ \left(p\cos(\theta),p\sin(\theta)\omega\right) : \theta\in (0,\theta_{u}],\omega\in S^{n-2} \}, $
$\partial_{2}U=\partial {\cal M} \cap \partial U$  and 
$\partial_{3}U= \{ \left(r\cos(\theta_{u}),r\sin(\theta_{u})\omega\right)\in \overline{\Omega} : r\in [0,p],\omega\in S^{n-2} \}, $
then $\partial U = \partial_{1}U \cup \partial_{2}U \cup \partial_{3}U,$  
$\tilde h_{2}\le 0\le f$  on $\partial_{1}U\setminus\{P\}$  and   $\tilde h_{2}<\tilde g^{[u]}<f$  on $\partial_{3}U$  (see (\ref{SPACESHIP}));  
then (\ref{Doggy}) and the general comparison principle imply 
\begin{equation}
\label{Jerry}
\tilde h_{2}<f \ \ \ \ \ {\rm in} \ \ \ \ \ U=\tilde A,
\end{equation}
where $A=\{r\left(\cos(\theta),\sin(\theta)\right)\in B_{2}({\cal O}_{2},p)\setminus \overline{M} : r\in (0,p), \theta\in (0,\theta_{u})\}$  
(see Figure \ref{FigureEleven}).  
Set ${\cal R}_{2} = \bigcup_{t=3\sigma/4}^{2\sigma/4}  l_{t}\left({\cal L}_{0}\right).$  
Now (\ref{STARSHINE})  implies $\tilde {\cal R}_{2}\subset U$  and so 
\begin{equation}
\label{Cab} 
f > \tilde h_{2}\ge -\frac{\beta\sigma}{4} \ \ \ \ \ {\rm on} \ \ \  {\cal R}_{2}. 
\end{equation}
Using (\ref{CatX}) and (\ref{Cab}), we see that if $a\in \left(p,\frac{2}{\lambda}-p\right)$  is close enough to $p,$  then 
$k_{2}(0,p)-k_{2}(0,a)<-\frac{\beta\sigma}{4}$  and therefore $f$  cannot be continuous at $P$  or at any point of $T=\{(0,p\omega)\in \Real^{n} : \omega\in S^{n-2} \}$.  
Notice that $f\in C^{0}\left(\overline{\Omega}\setminus T\right)$  (e.g. \cite{Lin}).  

\begin{figure}
\centering
\includegraphics{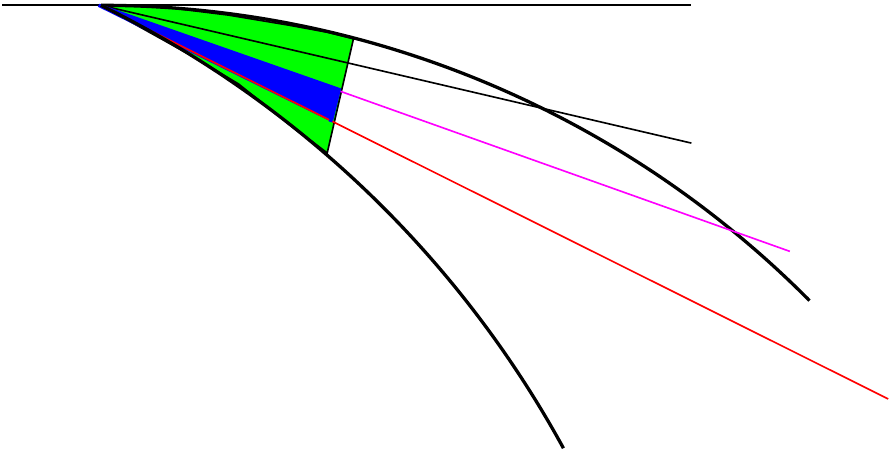}  
\caption{An illustration of ${\cal R}_{2}$ (blue region) and $A$ (green and blue regions)}
\label{FigureTwelve}
\end{figure}

\subsection{One singular point}
\label{OSP}

In this section, we will obtain a domain $\Omega$  and $\phi\in C^{\infty}(\Real^{n})$  such that $P\in\partial\Omega,$   the minimizer $f$  of 
(\ref{The_Functional}) is discontinuous at $P,$   $\partial\Omega\setminus \{P\}$  is smooth ($C^{\infty}$)
and $f\in C^{0}\left(\overline{\Omega}\setminus \{P\}\right).$  
This is accomplished by replacing ${\cal M}$  by  a convex set ${\cal G}$  such that $\partial{\cal G}\setminus \{P\}$  is smooth ($C^{\infty}$)  and 
${\cal G}\subset B_{n}\left({\cal O}_{n},p\right).$      We shall use the notation of \S \ref{CDO} throughout this section.  
We assume $p\in \left(0,\frac{1}{\lambda}\right)$  and set $P=(0,p,0,\dots,0).$    (We will no longer require Figure \ref{FigureEight} (b) to hold.) 

\begin{figure}
\centering
\includegraphics[width=5cm]{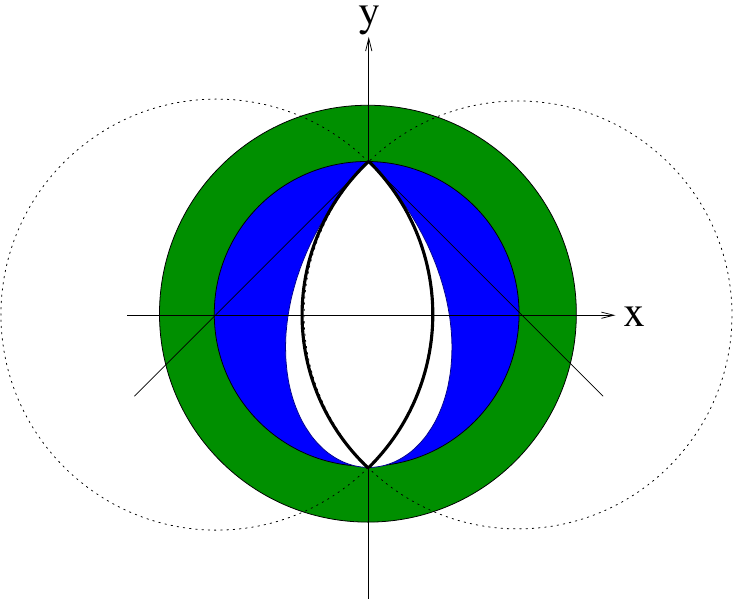}
\includegraphics[width=5cm]{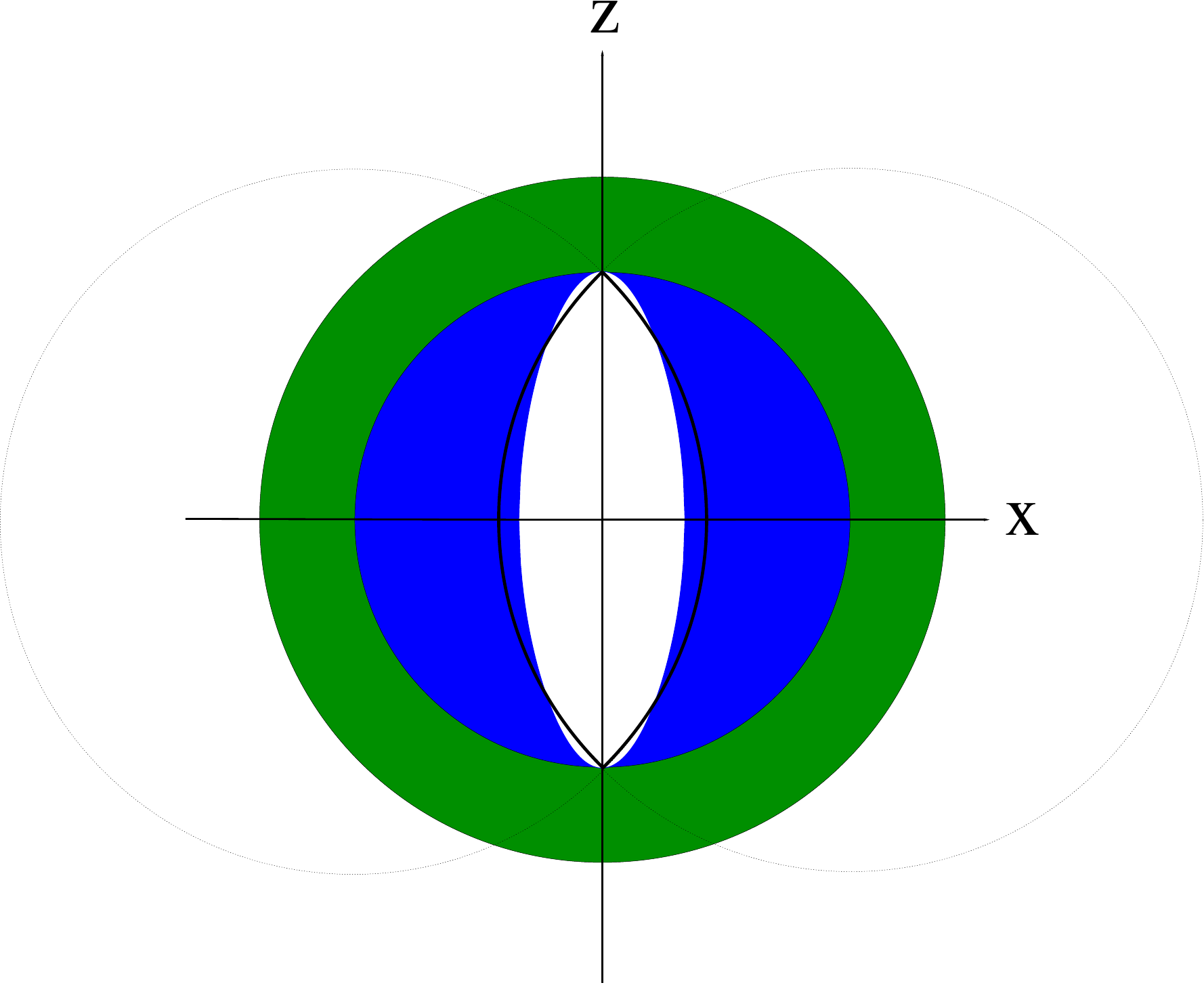}
\caption{(a) $\Pi_{1,2}\left(\Omega\right)$ \hspace{2cm} (b) $\Pi_{1,3}\left(\Omega\right)$ } 
\label{FigureThirteen}
\end{figure}

Let $\alpha>1,$  $n\ge 3,$  and $Y:\left[-\frac{\pi}{2\alpha},\frac{\pi}{2\alpha}\right] \times [0,\pi] \times S^{n-3}\to \Real^{n}$  be defined by 
\[
Y(\theta,\phi,\omega) = 2\cos(\alpha\theta)\sin(\phi)\left(\cos(\theta)\sin(\phi),\sin(\theta)\sin(\phi),\cos(\phi)\omega\right).
\]
Let $F:\Real^{n}\to \Real^{n}$  be given by $F\left(x_{1},\dots,x_{n}\right)=\left(\frac{x_{2}}{p},\frac{1-x_{1}}{p},\frac{x_{3}}{p},\dots,\frac{x_{n}}{p}\right)$
and define $X(\theta,\phi,\omega) = F\left(Y(\theta,\phi,\omega)\right)$  for $-\frac{\pi}{2\alpha}\le\theta\le \frac{\pi}{2\alpha}, \ 0\le\phi\le\pi, \ \omega\in S^{n-3}$
(see Figures \ref{FigureTwelve} and Figure \ref{FigureThirteen} with $n=3,$   $\alpha=2;$  the axes are labeled $x,y,z$  for $x_{1},x_{2},x_{3}$  respectively).
Let  ${\cal G}$   be the open, convex set whose boundary is the image of $X;$  that is, 
\[
\partial{\cal G}= \{X(\theta,\phi,\omega) : -\frac{\pi}{2\alpha}\le\theta\le \frac{\pi}{2\alpha}, \ 0\le\phi\le\pi, \ 
\omega\in S^{n-3} \}.
\]
Notice that $\partial{\cal G}\setminus \{P\}$  is a $C^{\infty}$  hypersurface in $\Real^{n}$  and $\partial {\cal G}\subset {\overline{B_{n}\left({\cal O}_{n},p\right)}}.$  
\begin{figure}
\centering
\includegraphics[width=4cm]{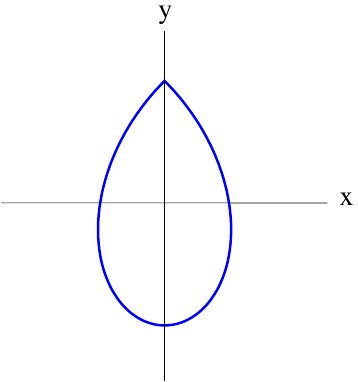}
\includegraphics[width=4cm]{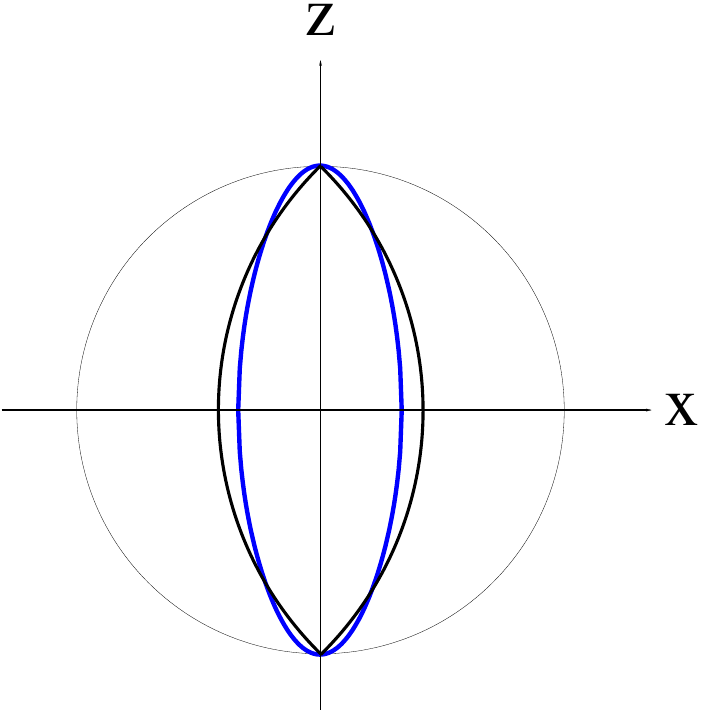}
\caption{$X\left(\theta,\frac{\pi}{2},1\right),\ \ X\left(\theta,\frac{1}{2}\arccos\right(1-\sec(\theta)\sec(2\theta)\left),1\right)$  }
\label{FigureTwelve}
\end{figure}

Let $\tau$  satisfy 
\[
0<\tau <  \min \left\{\frac{pr_{1}}{\sqrt{r_{2}^{2}-r_{1}^{2}}}, \frac{b(4p-b)}{4(2p-b)} \right\}.
\]
Set $\sigma=-\arctan(\tau/p)$  and $\alpha=\frac{\pi}{\pi+2\sigma}.$  
Then the tangent cones to $\partial {\cal G}$  and $\partial {\cal M}$  at $P$  are identical, $\cos(\sigma)>\frac{r_{1}}{r_{2}}$  and (\ref{Cat2}) holds 
for $u=(u_{1},u_{2})\in \partial B_{2}({\cal O}_{2},p)$   with $|u-w|<\delta_{1}.$   
If necessary by making $\tau>0$  smaller, we may assume $B_{n}\left({\cal O}_{n},\frac{a+p}{2}-b \right)\subset {\cal G}$  if $p<a<p+b.$

Now pick $a\in \left(p,\min\{p+b,1/\lambda\}\right)$  such that $k_{2}(0,p)-k_{2}(0,a)<-\frac{\beta\sigma}{4},$  as in (\ref{Cab}),  and define 
\begin{equation}
\label{BBBB}
\Omega=B_{n}\left({\cal O}_{n},a\right)\setminus \overline{\cal G}.
\end{equation}
Let 
\[
m>\max\{g_{1}(0,r_{1}), \frac{1}{2m_{0}}, \frac{\beta(2\sigma^{2}-\pi)}{8\sigma}\}
\]
and define $\phi\in C^{\infty}\left(\Real^{n}\right)$  so that $\phi=0$  on $\partial B_{n}\left({\cal O}_{n},a\right)$  and 
$\phi=m$  on $\partial {\cal G}$  and let $f$  be the variational solution of (\ref{ONE-A})-(\ref{ONE-B}).  
Notice that $f\in C^{2}(\Omega)$  satisfies (\ref{ONE-A})  and $f\in C^{0}\left(\overline{\Omega}\setminus \{P\}\right)$  (e.g. \cite{Lin}).

As in (\ref{Bingo}), let $B=B_{2}\left(\frac{p-r_{2}}{p}u,r_{2}\right).$   
Set $U_{0}=\{x\in\Omega : x\in \tilde B, x_{1}>0\}$  and   
$U=\{r\left(\cos(\theta),\sin(\theta)\omega\right)\in \Omega:  r\in (0,p), \theta\in (0,\theta_{u}),\omega\in S^{n-2} \}.$
Now $\tilde g^{[u]}=0$  on $\partial U_{0}\cap \partial \tilde B$  and $\tilde g^{[u]}\le g_{1}(0,r_{1})<m$  on 
$\partial U_{0}\cap \partial {\cal G}$  and so Lemma \ref{Four}, Lemma \ref{Mango} and (\ref{ONE-A}) imply $\tilde g^{[u]} \le f$  in $U_{0}$  
since $f$  minimizes the functional in (\ref{The_Functional}). 

As before, set $\partial_{1}U=\{ \left(p\cos(\theta),p\sin(\theta)\omega\right) : \theta\in [0,\theta_{u}],\omega\in S^{n-2} \},$
$\partial_{2}U=\partial {\cal G} \cap \partial U$  and 
$\partial_{3}U= \{ \left(r\cos(\theta_{u}),r\sin(\theta_{u})\omega\right)\in \overline{\Omega} : r\in [0,p],\omega\in S^{n-2} \}.$
Then $f\ge 0$  on $\partial_{1}U\setminus\{P\},$  $\partial U = \partial_{1}U \cup \partial_{2}U \cup \partial_{3}U,$  
$\tilde h_{2}\le 0\le f$  on $\partial_{1}U,$   $\tilde h_{2}< m=\phi$  on $\partial_{2}U$  and  $\tilde h_{2}<\tilde g^{[u]}<f$  
on $\partial_{3}U;$  Lemma \ref{Four} implies that (\ref{Cab}) continues to hold.  
Then (\ref{CatX}) and (\ref{Cab}) imply $f$  is discontinuous at $P$  since $k_{2}(0,p)-k_{2}(0,a)<-\frac{\beta\sigma}{4}.$

\section{The Concus-Finn conjecture}
\label{CFcondition}
For the moment, assume $n=2.$
In approximately 1970, Paul Concus and Robert Finn conjectured that if $\kappa\ge 0,$  $\Omega\subset \Real^{2}$  has a corner at $P\in\partial\Omega$  
of (angular) size $2\alpha,$   $\alpha\in\left(0,\frac{\pi}{2}\right),$   $\gamma:\partial\Omega\setminus\{P\}\to [0,\pi]$  and 
$|\frac{\pi}{2}-\gamma_{0}|>\alpha,$   where 
\begin{equation}
\label{FROZEN}
\lim_{\partial\Omega\ni x\to P} \gamma(x)=\gamma_{0},
\end{equation}
then a function $f\in C^{2}(\Omega)\cap C^{1}\left(\overline{\Omega}\setminus \{P\}\right)$  which satisfies 
\begin{eqnarray}
\label{TEN-A}
{\rm div}\left( Tf\right) &  = &  \kappa f \ \ \ \ \ {\rm in} \ \ \Omega, \\ 
Tf\cdot\eta & = & \cos(\gamma) \ \ \ \ \ {\rm  on} \ \ \partial\Omega\setminus \{P\}, 
\label{Fletch}
\end{eqnarray}
must be discontinuous at $P;$   here $\eta(x)$  is the exterior unit normal to $\Omega$  at $x\in\partial\Omega\setminus \{P\}.$
A generalization (including the replacement of (\ref{TEN-A}) by (\ref{ONE-A}))  of this conjecture in the case $\gamma_{0}\in \left(0,\pi\right)$  
was proven in \cite{CFC}.

In the situation above with $\alpha\in\left(\frac{\pi}{2},\pi\right),$  the ``nonconvex Concus-Finn conjecture'' states that if 
$|\frac{\pi}{2}-\gamma_{0}|>\pi-\alpha,$  then the capillary surface $f$  with contact angle $\gamma$   must be discontinuous at $P.$
A generalization (including the replacement of (\ref{TEN-A}) by (\ref{ONE-A})) of this extension of the Concus-Finn conjecture  
in the case $\gamma_{0}\in \left(0,\pi\right)$  was proven in \cite{NCFC}.
Both \cite{CFC} and \cite{NCFC} include the possibility of differing limiting contact angles; that is, the following limits 
\[
\lim_{\partial^{+}\Omega\ni x\to P} \gamma(x)=\gamma_{1} \ \ \ \ \ {\rm and} \ \ \ \ \ \lim_{\partial^{-}\Omega\ni x\to P} \gamma(x)=\gamma_{2}
\]
exist, $\gamma_{1}, \gamma_{2}\in (0,\pi)$   and $\gamma_{1}\neq \gamma_{2}.$  
Here $\partial^{+}\Omega$  and $\partial^{-}\Omega$  are the two components of $\partial\Omega\setminus \{P,Q\},$
where $Q\in \partial\Omega\setminus \{P\}.$  When $\gamma_{1}\neq \gamma_{2},$  the necessary and sufficient (when $\alpha\le \frac{\pi}{2}$)  or necessary 
(when $\alpha>\frac{\pi}{2}$)  conditions for the continuity of $f$  at $P$  become slightly more complicated. 

The cases where $\gamma_{0}=0,$  $\gamma_{0}=\pi,$  $\min\{\gamma_{1},\gamma_{2}\}=0$  and $\max\{\gamma_{1},\gamma_{2}\}=\pi$  remain unresolved.  
If we suppose for a moment that the nonconvex Concus-Finn conjecture with limiting contact angles of zero or $\pi$  is proven, then the discontinuity 
of $f$  at $P$  in \S \ref{BLUE} follows immediately from the fact that $f<\phi$  in a neighborhood  in $\partial\Omega\setminus \{P\}$  of $P$ 
since then Lemma \ref{Three} implies $\gamma_{0}=0$  and therefore $|\frac{\pi}{2}-\gamma_{0}|>\pi-\alpha.$
In this situation (i.e. the solution $f$  of a Dirichlet problem satisfies a zero (or $\pi$)  contact angle boundary condition near $P$), 
establishing the discontinuity of $f$  at $P$  would be much easier and a much larger class of domains $\Omega$  with a nonconvex corner 
(i.e. $\alpha>\frac{\pi}{2}$)  at $P$  would have this property.
For example, if $\Omega$  is a bounded locally Lipschitz domain in $\Real^{2}$  for which (\ref{PIZZA}) holds, $f\in C^{2}(\Omega)$  is a generalized solution of  
(\ref{ONE-A})-(\ref{ONE-B}) (and $H$  need not vanish) and $\phi$  is large enough near $P$  (depending on $H$  and the maximum of $\phi$  outside 
some neighborhood of $P$) that $f<\phi$  on $\partial\Omega\setminus \{P\}$  near $P,$   
then the fact that $\gamma_{0}=0$  (Lemma \ref{Three})  together with the nonconvex Concus-Finn conjecture would imply that $f$  is discontinuous at $P.$

Now consider $n\in\Natural$  with $n\ge 3.$  
Formulating generalizations of the Concus-Finn conjecture in the ``convex corner case'' (i.e. 
$\Omega\cap B_{n}(P,r)\subset \{X\in\Real^{n}: (X-P)\cdot \mu>0\}$  for some $\mu\in S^{n-1},$  $P\in\partial\Omega$  and $r>0$) and in other 
cases where $\partial\Omega$  is not smooth at a point $P\in\partial\Omega$  may be complicated because the geometry of 
$\partial\Omega\setminus \{P\}$  is much more interesting when $n>2.$  
Establishing the validity of a generalization of the Concus-Finn conjecture for solutions of (\ref{ONE-A}) \& (\ref{Fletch}) 
when $n>2$  is probably significantly harder than doing so when $n=2.$  

Suppose we knew that a solution $f$  of (\ref{ONE-A}) \& (\ref{Fletch}) is necessarily discontinuous at a ``nonconvex corner'' $P\in\partial\Omega$  
when $\gamma_{0}=0,$  where $\gamma_{0}$  is given by (\ref{FROZEN}).  
In this case, a necessary condition for the continuity of $f$  at $P$  would be that $\limsup_{\partial\Omega\ni X\to P} Tf\left(X\right)\cdot \eta(X) >0$  
and $\liminf_{\partial\Omega\ni X\to P} Tf\left(X\right)\cdot \eta(X)<\pi.$
Then the arguments in \S \ref{Happy} could be made more easily and the conclusion 
that $f$  is discontinuous at $P$  would hold in a much larger class of domains $\Omega;$  here, of course, we use the ridge point $P$ in \S \ref{Happy} 
as an example of a ``nonconvex corner'' of a domain in $\Real^{n}.$  
The primary difficulty in proving in \S \ref{Happy} that $f$  is discontinuous at $P$  is establishing (\ref{Cab}); a more ``natural'' generalization 
of $\Omega\subset\Real^{2}$  in \S \ref{BLUE} would be 
\[
\Omega^{*} = \{(x\omega_{1},y,\omega_{2},\dots,\omega_{n-1})\in \Real^{n} : (x,y)\in B_{2}\left({\cal O}_{2},a\right)\setminus \overline{M}, \ 
\omega\in S^{n-1}\}.
\]
However, the use of Lemma \ref{Mango} to help establish (\ref{Cab}) in $\Omega^{*}$  is highly problematic.  
On the other hand, an n-dimensional ``Concus-Finn theorem'' for a nonconvex conical point (e.g. $P\in\partial\Omega^{*}$) would only require an 
inequality like (\ref{Taxi}) to prove that $f<\phi$  on $\partial\Omega\setminus\{P\}$   near $P$  and hence that $f$  is discontinuous at $P;$ 
the replacement of (\ref{AAAA})  by (\ref{BBBB}) in order to obtain a $\Omega$  such that $\partial\Omega\setminus\{P\}$  is $C^{\infty}$  would be 
unnecessary.


\end{document}